\documentclass[12pt]{article}

\usepackage{amsfonts}

\newtheorem{thm}{Theorem}[section]
\newtheorem{lemma}[thm]{Lemma}
\newtheorem{cor}[thm]{Corollary}

\newtheorem{prop}[thm]{Proposition}

\newcommand{\be}{\begin{equation}}
\newcommand{\ee}{\end{equation}}
\newcommand{\openbox}{\leavevmode
  \hbox to8pt{\hfil\vrule\vbox to6pt{\hrule width6pt\vfil\hrule}\vrule}}

\newcommand{\qed}{\hbox to5pt{ } \hfill \openbox\bigskip\medskip}

\newcommand{\Fp}{\mathbb F _p}

\newcommand{\Sf}{\mathbb S}
\newcommand{\ve}[1]{\mathbf{#1}}

\newcommand{\cF}{\mbox{$\cal F$}}

\newcommand{\cH}{\mbox{$\cal H$}}

\newcommand{\cN}{\mbox{$\cal N$}}

\newcommand{\cV}{\mbox{$\cal V$}}

\newcommand{\cS}{\mathbb S}

\newcommand{\R}{\mathbb R}

\newcommand{\F}{\mathbb F}

\title{Equality in the linear algebra bound}
\author{G\'abor Heged\"{u}s
\\{\normalsize  \'Obuda University}
\\{\normalsize B\'ecsi \'ut 96/B, Budapest, Hungary, H-1032}
\\{\normalsize {\tt hegedus.gabor@uni-obuda.hu}}
\\ Lajos R\'onyai
\\ {\normalsize HUN-REN Institute for Computer Science and Control}
\\ {\normalsize and }
\\ {\normalsize Dept. of Algebra and Geometry,}
\\ { \normalsize Budapest University of Technology and Economics}
\\ {\normalsize \tt lajos@ilab.sztaki.hu }
}

\begin{document}
\maketitle

\begin{abstract}
We study some examples when there is actually an equality in the linear
algebra bound. When the vectors considered span in fact the entire space. 
We would like to point out that in some cases this provides some interesting 
extra information about the extremal configuration. We obtain results on set
families satisfying conditions on pairwise intersections, or Hamming
distances. Also, we have an application to 2-distance sets in Euclidean 
spaces. 
\end{abstract}
\medskip

\section{Introduction}

The linear algebra bound method in discrete mathematics works by assigning
vectors
to some objects of interest in a linear space $V$ over a field. By
proving
that the vectors are linearly independent, we obtain that the number of
objects is at most the dimension of $V$. Here we focus on the case when the
vectors in fact form a basis of $V$. In some situations this basis property
allows us to obtain interesting additional information on the extremal
configurations. In this paper we intend to give some examples of this 
phenomenon.
\medskip

Throughout the paper $n$ and $q>1$ are positive integers and write 
$[n]:=\{1,2,,\ldots, n\}$, and $[0,q-1]=\{0,1,\ldots, q-1\}$.
Let $\cH\subseteq [0,q-1]^n$ and 
let $\ve h_1,\ve h_2\in \cH$ be two elements of the vector system $\cH$.  
Let $d_H(\ve h_1,\ve h_2)$ stand for the Hamming distance of the
vectors $\ve h_1,\ve h_2\in \cH$:
$$
d_H(\ve h_1,\ve h_2):=|\{i\in [n]:~  (\ve h_1)_i\neq (\ve h_2)_i\}|.
$$

Denote by $D(\cH)$ the following set of Hamming distances: 
$$
D(\cH):=\{ d_H(\ve h_1,\ve h_2):~ \ve h_1,\ve h_2\in \cH, \ve h_1\neq \ve h_2 \}.
$$  

Delsarte proved the following well-known upper bound for the size of 
a vector system with $s$ distinct Hamming distances (see in (1.2) of
\cite{D3}).

\begin{thm} \label{Delsarte}
Let $0<s\leq n$, $q>1$ be positive integers, 
and $$L=\{\ell_1 ,\ldots ,\ell_s\}\subseteq [n]$$
be a set of $s$ 
positive integers. Suppose that $\cH\subseteq [0,q-1]^n$ and  
that $d_H(\ve h_1,\ve h_2)$ is in $L$ for each pair of distinct 
vectors $\ve h_1,\ve h_2\in \cH$. Then
$$
|\cH|\leq \sum_{i=0}^s {n\choose i} (q-1)^i.
$$
\end{thm}


In \cite{H} the first author studied the maximal size of families with a 
unique Hamming distance  between distinct members of the family: 

\begin{thm} \label{old}
Let $\cF=\{F_1, \ldots , F_m\}$ be  a family of subsets of $[n]$ such that 
there exists a positive integer $\lambda>0$ with 
$d_H(\ve v_i,\ve v_j)=\lambda$ for 
each $i\neq j$. Here $\ve v_i$
is the characteristic vector of the set $F_i\in \cF$. 
Suppose further that $\lambda\neq \frac{n+1}{2}$.
Then $m=|\cF|\leq n$.
\end{thm}

Please note that the condition on $\lambda$ can not be dropped. Let 
${\cal D} \subset 2^{[4v-1]}$
be a $(4v-1,2v-1, v-1)$ Hadamard design for some positive integer $v$ and
set $\cF:= {\cal D}\cup\{[4v-1]\}$. Then $d_H(\ve v_i,\ve v_j)=2v$ whenever 
$\ve v_i$ and $\ve v_j$ represent different sets from 
$\cF$, and $|\cF |=4v=n+1$.     

Our first result generalizes to vector systems and 
gives a modular 
version of Theorem \ref{old} as well as of Theorem 6 in \cite{HHY}. 
In the proof we use the basis property of the vectors involved in the linear 
algebra bound.

\begin{thm} \label{main}
Let $n$ and $q>1$ be positive integers. Let $p\geq q$, $p$ be a prime, 
and  assume, that $n\not\equiv 0 \pmod p$. 
Let $\cH\subseteq [0,q-1]^n\subseteq {\Fp}^n$ 
be a vector system such that there exists a positive integer 
$\lambda\not\equiv 0 \pmod p$ with $d_H(\ve h,\ve g)\equiv \lambda \pmod p$ 
for each distinct  pair $\ve h,\ve g\in \cH$.  
Finally suppose that $q\lambda \not\equiv n(q-1)+1 \pmod p$.
Then $|\cH|\leq n(q-1)$.
\end{thm}

\medskip


\medskip


Let $\langle \ve x,\ve y\rangle $ stand for  the standard scalar product 
on $\R^n$.  Let $A(\cV)$ denote the set of scalar products between 
different vectors of $\cV\subseteq \R^n$:
$$
A(\cV):=\{\langle \ve p_1,\ve p_2\rangle :~ \ve p_1,\ve p_2\in \cV,~ \ve p_1\ne \ve p_2\}.
$$

A {\em spherical $s$-distance set} is a 
subset $\mbox{$\cal V$}\subseteq {\Sf}^{n-1}$ such 
that $|A(\mbox{$\cal V$})|\leq s$.
Let $n,s\geq 1$ be integers. Set
$$
M(n,s):={n+s-1\choose s}+{n+s-2\choose s-1}.
$$
Delsarte, Goethals and  Seidel investigated spherical $s$-distance sets. 
They  proved a general upper bound 
in \cite{DGS}.

\begin{thm} \label{DGGupper}
Suppose that $\cV\subseteq {\Sf}^{n-1}$ is a spherical set 
satisfying $|A(\cV)|\leq s$.
Then
$$
|\cV|\leq M(n,s).
$$
\end{thm}

We point out the following special case of Theorem \ref{DGGupper}.
\begin{cor} \label{2dist}
Suppose that $\cV\subseteq {\cS}^{n-1}$ is a set satisfying 
$|A(\cV)|\leq 2$. Then
$$
|\cV|\leq \frac{n(n+3)}{2}.
$$
\end{cor}

Musin proved a related result in Theorem 1 of \cite{M}. 
He employed  the linear algebra bound.
The parameters of large spherical two-distance set are constrained.
For instance, Neumaier proved the following result in Corollary 5 of 
\cite{N}.

\begin{thm} \label{Neu}
Let $\cV$ be an $n$-dimensional two-distance set with distances 
$d_1,d_2$ (where $d_1<d_2$). If $|\cV|>\max(2n+1,5)$, then 
there exists an integer $m$ such 
that $\frac{d_1^2}{d_2^2}=\frac{(m-1)}{m}$.
\end{thm}

In a similar spirit, with the linear algebra bound method we
exhibit here an algebraic relation of the parameters of a maximal
spherical 2-distance set:  

\begin{thm} \label{main2}
Let $n>1$ be a  positive integer and 
define $N:=\frac{n(n+3)}{2}$.  
Let $\cV:=\{\ve v_1, \ldots , \ve v_N\}\subseteq \cS^{n-1}$ be 
a set of unit vectors such that $A(\cV)=\{a,b\}$, where 
$a\neq 1$,  $b\neq 1$.
 Then 

\begin{equation} \label{algebrai}
N\Big(ab + \frac 1n\Big)=(1-a)(1-b).
\end{equation}
\end{thm}


Please note that Barg et al. in Theorem 2.4 of \cite{BGOY} obtained the same
formula (\ref{algebrai}), as one of two alternative equations connecting 
$a$, $b$, $n$, and
$N$, in the related setting of two distance unit norm tight frames. Our
result is about maximal spherical two-distance sets.


As our next example of the use of the basis property, 
we consider a modulo $p$-uniform 
set family with a modular intersection condition. We obtain a modular 
variant of the equation valid for symmetric designs: 

\begin{thm} \label{main4}
Let $p$ be a prime and $k,\lambda$ be non-negative integers, 
$n,k,k-\lambda$ are not
divisible by $p$.   
Let $\cF=\{F_1, \ldots , F_n\}$ be  a family of subsets
of $[n]$ such that $|F_i|\equiv k \pmod p$ for every $i$, and  
$|F_i\cap F_j|\equiv \lambda \pmod p$ for each $i\neq j$. 
Then we have $k(k-1)\equiv \lambda (n-1) \pmod p$. Moreover for the degree
$d_i$ of every $i\in [n]$ with respect to $\cF$ we have 
$d_i\equiv k \pmod p$. 
\end{thm}

\noindent
{\bf Remarks.}
1. Ryser designs provide infinitely many examples when 
Theorem \ref{main4} applies, and $\cF$ is not a uniform (but 
of course, $p$-uniform)
family.  Let $p>2$ be a prime and $r$ be a prime of the form $r=dp+1$. By
Dirichlet's theorem on primes in arithmetic progressions there are 
infinitely many such primes $r$. Consider 
an $(r^2+r+1,r+1,1)$ projective plane and let $\cF$ be the type 1 
$\lambda$-design obtained from it. Then $\cF$ has $n=r^2+r+1$ points and 
$n$ blocks, the block sizes are $2r$ and $r+1$, and $\lambda=r$. The
residues modulo $p$ of $n,k,\lambda$, are 3,2,1, respectively. \\
2. For type 1 $\lambda$-designs the congruence of Theorem \ref{main4}
follows very easily. Indeed, let we have a symmetric $(n,k',\lambda')$
design. We have then 
\begin{equation} \label{egyenlet}
\lambda'(n-1)=k'(k'-1).  
\end{equation}
The resulting $\lambda$-design will have $n$ points, block sizes
$2(k'-\lambda')$, and $k'$, and $\lambda=k'-\lambda'$. Also, $p$ is a prime
divisor of $k'-2\lambda'$. 
The congruence to be verified is 
$$ (k'-\lambda')(n-1)=k'(k'-1) \pmod p ,$$
which follows at once from (\ref{egyenlet}) and 
$k'-\lambda'\equiv \lambda' \pmod p$. 

\medskip

Finally, with the polynomial method  we give another proof of a result of 
Ryser, the extremal case
of Theorem 1.1 from \cite{Ryser}. Our 
argument derives perhaps more directly essentially the same algebraic facts 
as Ryser's proof. 

\begin{thm} (Ryser)\label{main5}
Let $\lambda$ denote a positive integer,
 $\cF=\{F_1, \ldots , F_n\}$ a family of subsets
of $[n]$, with $|F_i\cap F_j|=\lambda$ for each $i\neq j$. 
Suppose further that $|F_i|>\lambda$ for each $i$. Then one of the
following statements holds: \\
(A) There is a positive integer $r$ such that any point of $[n]$ is
contained in exactly $r$ elements of $\cF$ and $|F|=r$ for any
$F\in \cF$; \\
(B) There exist different positive integers $r, r'$ such that $r+r'=n+1$ and 
any point of $[n]$ occurs in either $r$ or $r'$ elements of $\cF$.
\end{thm}

In Section 2 we prove  our results. In Section 3 a problem for further
research is outlined.  

\section{Proofs}

The proofs are based on the linear algebra bound method. 
In particular, we shall use the following Determinant Criterion (see e.g.
Proposition 2.7 in \cite{BF}). 
\begin{prop} \label{det} (Determinant Criterion)
Let $\F$ denote an arbitrary field. Let $f_i:\Omega \to \F$ be functions 
and $\ve v_j\in \Omega$ elements for each $1\leq i,j\leq m$  such that 
the $m \times m$ matrix $B=(f_i(\ve v_j))_{i,j=1}^m$
is non-singular. Then $f_1,\ldots ,f_m$ are linearly independent functions 
of the space $\F^{\Omega}$. 
\end{prop}
It is easy to verify the following Lemma.

\begin{lemma} \label{meret}
Let $p\geq q$, $p$ be a prime, $\ve f\in [0,q-1]^n\subseteq {\Fp}^n$ 
be an arbitrary 
vector. Let $\ve j=(j,\ldots, j) \in [0,q-1]^n$ denote a constant vector. 
Then in $\F_p$ we have 
$$
\sum_{j=0}^{q-1} d_H(\ve f, \ve j)=n(q-1).
$$
\end{lemma}



\medskip

{\bf Proof of Theorem \ref{main}:}
Let $a\in [0,q-1]$, and  $l_a(x)\in\Fp[x]$ be the univariate polynomial with 
minimal degree such that $l_a(a)=0$ and $l_a(b)=1$ for 
each $b\in [0,q-1]$, $b\neq a$. Clearly $\deg(l_a)\leq q-1$. 

Let $\ve a=(a_1, \ldots ,a_n)\in \cH$ be a vector and define the 
multivariate polynomial
$$
f_{\ve a}(x_1, \ldots ,x_n):= \sum_{i=1}^n l_{a_i}(x_i)-\lambda\in \Fp[x_1, \ldots ,x_n].
$$

Then it is easy to verify using Proposition \ref{det}  and the condition 
$\lambda\not\equiv 0 \pmod p $ that the set of 
polynomials $\{f_{\ve a}:~ \ve a\in \cH\}$ is linearly independent 
over $\F_p$ (see also the proof of Theorem 5 in \cite{BSW}).
Also,  each polynomial $f_{\ve a}$ has the property that the monomials 
involved in $f_{\ve a}$ depend on at most one indeterminate $x_i$. 
Consequently  all the polynomials $f_{\ve a}$, 
where $\ve a\in \cH$, appear in the linear span of $n(q-1)+1$ monomials 
which depend on at most one indeterminate $x_j$ and the 
exponent of $x_j$ is at most $q-1$. Here we used that $\deg(l_a)\leq q-1$. 

Suppose that $|\cH|=n(q-1)+1$. Then  $\{f_{\ve a}:~ \ve a\in \cH\}$ is a 
basis of the $\F_p$-linear space generated by the monomials  
$$\{1\}\cup \{x_i^j:~ 1\leq i\leq n,~ 1\leq j\leq q-1\}.$$ 
This means that 
we can write up the constant  polynomial  $1\in  \Fp[\ve x ]$ as a linear 
combination of the polynomials  from $\{f_{\ve a}:~ \ve a\in \cH\}$:
\begin{equation}  \label{comb}
1=\sum_{\ve a\in\cH} \alpha_{\ve a}f_{\ve a},
\end{equation}
where $ \alpha_{\ve a}\in \Fp$. 

Substituting $\ve b\in \cH$ into equation (\ref{comb}), we see that
$$ 
\alpha_{\ve b}=- \frac {1}{\lambda}
$$ 
for each $\ve b\in \cH$. It follows that in $\F_p[\ve x]$ we have 
\begin{equation}  \label{comb2}
1=- \sum_{\ve a\in\cH} \frac {1}{\lambda}f_{\ve a}.
\end{equation}

Now substitute the vector $\ve j:=(j, \ldots ,j)\in \Fp^n$ for 
each $j\in [0,q-1]$ into equation (\ref{comb2}). This gives us
\begin{equation}  \label{comb3}
-\lambda= \sum_{\ve a\in\cH} f_{\ve a}(\ve j) \pmod p.
\end{equation}
Observe that 
$$
f_{\ve a}(\ve j)=d_H(\ve a,\ve j)-\lambda \pmod p
$$
for each $j\in [0,q-1]$. Adding up equations   (\ref{comb3}), we obtain
$$
-\lambda\cdot q =\sum_{j\in [0,q-1]}\Big(\sum_{\ve a\in\cH}  
f_{\ve a}(\ve j)\Big)=
$$
$$
=\sum_{\ve a\in\cH}  \Big(\sum_{j\in [0,q-1]} 
f_{\ve a}(\ve j)\Big)=\sum_{\ve a\in\cH}  
\sum_{j\in [0,q-1]} \Big( d_H(\ve a,\ve j)-\lambda \Big) \pmod p.
$$
Now it follows from Lemma \ref{meret} that
$$
\sum_{j\in [0,q-1]} \Big(d_H(\ve a,\ve j)-\lambda\Big)
=n(q-1)-\lambda\cdot q \pmod p, 
$$
and therefore 
$$
-\lambda\cdot q=\sum_{\ve a\in\cH} \Big( n(q-1)-\lambda\cdot q \Big)
=|\cH| (n(q-1)-\lambda\cdot q) 
\pmod p.
$$
As $|\cH|=n(q-1)+1$, finally we obtain
$$
-\lambda\cdot q=(n(q-1)+1)\cdot(n(q-1)-\lambda\cdot q) \pmod p
$$
It easy to verify from this equation that 
$\lambda\cdot q= n(q-1)+1 \pmod p$, which is in contradiction with 
the last assumption of the theorem. One uses here that 
$n (q-1) \not \equiv 0 \pmod p$. \qed

{\bf Proof of Theorem \ref{main2}:}
Consider the set $\cN(n)$ of monomials in variables $x_1,\ldots ,x_n$
which have degree at most 2 and have degree at most 1 in $x_1$.
Then clearly we have $1\in \cN(n)$, and
$$
|\cN(n)|=\frac{n(n+3)}{2}.
$$

We set $N:=\frac{n(n+3)}{2}$. 
Let $\cV:=\{\ve v_1, \ldots , \ve v_N\}\subseteq \cS^{n-1}$ be a 
collection of unit vectors such that $A(\cV)=\{a,b\}$, where 
$a\neq 1$,  $b\neq 1$. 
Consider the real polynomial
$$
g(x_1,\ldots ,x_n)=(\sum_{i=1}^n x_i^2) -1\in \R[x_1, \ldots ,x_n].
$$
Define the multivariate polynomial
$$
P_{m}(\ve x):=(\langle \ve x, 
\ve v_m\rangle -a)\cdot (\langle \ve x, \ve v_m\rangle -b)\in \R[\ve x],
$$
for each $1\leq m\leq N$, where $\langle \ve x, \ve y\rangle $ denotes the 
standard inner product on $\R^n$. 

We note first that any     $\ve s\in {\Sf}^{n-1}$ is a zero of the equation 
\begin{equation}  \label{sphere} 
x_1^2=1-\sum_{i=2}^n x_i^2.
\end{equation}
Let $Q_m$ denote the polynomial obtained by writing  $P_m$ as a linear 
combination of monomials and replacing each occurrence of $x_1^2$ 
by a linear combination of other 
monomials, using the relation (\ref{sphere}).

We have  $g(\ve s)=0$ for each $\ve s\in {\Sf}^{n-1}$, 
hence $Q_m(\ve s)=P_m(\ve s)$ also holds. 

It is easy to verify from  the Determinant Criterion, with the choices of 
$\F:=\R$, $\Omega :={\Sf}^{n-1}$ and $f_m:=Q_m$ for each $m$, that the set of 
polynomials $\{Q_m:~ 1\leq m\leq N\}$ is linearly independent.

Then it is easy to check that  we can write $Q_m$ as a linear combination of 
monomials in the form
$$
Q_m(\ve x)=\sum d_{\alpha}x^{\alpha},
$$ 
where $d_{\alpha}\in \R$ are real coefficients, and 
$x^{\alpha}:=x_{1}^{\alpha_{1}}\cdot \ldots \cdot x_{n}^{\alpha_{n}}\in
\cN(n)$.   
This follows immediately from relation (\ref{sphere}). 

Let $V$ denote the  vector space, which is generated  by the set of 
monomials $\cN(n)$.
Since $\{Q_m:~ 1\leq m\leq N\}$ is a set of  linearly independent  
polynomials in the vector space $V$  and $N=|\cN(n)|$, 
the set $\{Q_m:~ 1\leq m\leq N\}$  is a basis in the linear space $V$. This 
implies that we can write the constant  polynomial  
$1$ as a linear combination of the polynomials  $\{Q_m:~ 1\leq m\leq N\}$:
\begin{equation}  \label{comb22}
1=\sum_{m=1}^N \alpha_{m}Q_m,
\end{equation}
where $ \alpha_{m}\in \R$ for each $m$. 

Then substituting vector $\ve v_m\in \cV$ into equation 
(\ref{comb22}) we obtain that
$ \alpha_{m}=\frac {1}{(1-a)(1-b)}$ for each $1\leq m\leq  N$. We have the 
following equation  
\begin{equation}  \label{comb222}
1=\sum_{m=1}^N  \frac {1}{(1-a)(1-b)} Q_m, 
\end{equation}
as equality of polynomials. Using the fact that $Q_m(\ve s)=P_m (\ve s) $
whenever $\ve s \in  {\Sf}^{n-1}$, we have 
\begin{equation}  \label{P_i}
(1-a)(1-b)=\sum_{m=1}^N P_m(\ve s)
\end{equation}
for $\ve s \in  {\Sf}^{n-1}$.

Let $\ve e_i  \in \cS^{n-1}$, $1\leq i\leq n$ be the standard 
basis vectors: the $i$-th coordinate of $\ve e_i$ is 1, the others are 0. 
First we observe that

\begin{equation} \label{Pe}
P_m(\pm \ve e_i))=((\ve v_m)_i)^2 \mp (a+b)
                      (\ve v_m)_i+ab.   
\end{equation}

\noindent
We simplify the relation  
\begin{equation}
\sum _{m=1}^N P_m(\ve e_i)= \sum _{m=1}^N P_m(-\ve e_i)
\end{equation}
which holds because both sides are $(1-a)(1-b)$.  
We obtain  for every 
$1\leq i\leq n$  
\begin{equation}\label{e_i}
(a+b)\sum _{m=1}^N (\ve v_m) _i=0. 
\end{equation}

Next, using (\ref{P_i}), (\ref{Pe}), and (\ref{e_i}) for every $i$, 
$1\leq i\leq n$ we have 
\begin{equation}
(1-a)(1-b)=\sum _{m=1}^N P_m(\ve e_i)=\sum _{m=1}^N (((\ve v_m)_i)^2 - (a+b)
                      (\ve v_m)_i+ab)=
\end{equation}  
\begin{equation}
=\sum _{m=1}^N((\ve v_m)_i)^2  +Nab, 
\end{equation}
and hence 
\begin{equation}
\sum _{m=1}^N((\ve v_m)_i)^2=(1-a)(1-b)-Nab.
\end{equation}
Now add these up for $i=1,\ldots, n$:
\begin{equation} \label{N}
N=\sum_{i=1}^n\sum _{m=1}^N((\ve v_m)_i)^2= n(1-a)(1-b)-nNab, 
\end{equation}
where the first equality holds because the $\ve v _m$ are unit vectors in
$\R^n$: 
$$
\sum_{i=1}^n ((\ve v_m)_i)^2 =1.
$$
From (\ref{N}) dividing by $n$ and rearranging gives the desired equation 
$$
(1-a)(1-b)=N(\frac 1n +ab).
$$
\qed

{\bf Proof of Theorem \ref{main4}:}
Let $\ve v_i  \in \F_p ^n $, $1\leq i\leq n  $ denote the characteristic 
vector of $F_i$. 
Consider the polynomials 
$$
g_{i}(\ve x):=\langle \ve x,
\ve v_i\rangle -\lambda  \in \F_p [x_1,\ldots, x_n].
$$
Let $d_i$ be the degree of $i\in [n]$ with respect to the set family
$\cF$. Write $d=d_1$. Without loss of generality we may assume that 
$1\in F_1,\ldots , F_d$.

For $1\leq i\leq d$ consider the polynomial $h_i(\ve x)$ obtained from 
$g_i(\ve x)$ by substituting 
$-x_2-...-x_n+k $
for the variable $x_1$. We obtain that the polynomials 
\begin{equation} \label{independent}
h_1,\ldots ,h_d, g_{d+1},\ldots, g_n 
\end{equation}
are in $ \F_p[x_2,\ldots ,x_n]$.
We easily see that $h_i(\ve v_i)\equiv g_j(\ve v_j)\equiv k-\lambda \pmod
p$ for $1\leq i\leq d$, $d<j\leq n$, 
and $h_i(\ve v_j)\equiv g_i(\ve v_j)\equiv 0 \pmod p$ if $i\not= j$. 
These congruences imply that the polynomials in (\ref{independent}) are 
linearly 
independent over $\F_p$, and hence form a basis of the subspace of
polynomials of degree at most 1 in $\F_p[x_2,\ldots, x_n]$. As a 
consequence, 
the constant polynomial $k-\lambda$ is a linear combination 
\begin{equation} 
 k-\lambda=\alpha_1h_1(\ve x)+\cdots +\alpha_d h_d(\ve x)+
\alpha_{d+1}g_{d+1}(\ve x)+\cdots
+\alpha_n g_n(\ve x).
\end{equation}
By substituting $\ve v_i$ into the above equation, we obtain that 
$\alpha_i=1$ in $\F_p$ for every index $i$:  

\begin{equation} \label{pontos}
 k-\lambda=h_1(\ve x)+\cdots +h_d(\ve x)+ g_{d+1}(\ve x)+\cdots
+ g_n(\ve x).
\end{equation}
Comparing here the constant terms gives 
\begin{equation} \label{congr}
 k-\lambda\equiv -\lambda n + kd \pmod p.
\end{equation}
Essentially the same argument gives similar 
congruences for $1\leq i\leq
n$:
\begin{equation}
 k-\lambda\equiv -\lambda n + kd_i \pmod p.
\end{equation}
Using the assumption $k\not\equiv 0 \pmod p$ we obtain that 
$d_i\equiv d \pmod p $ for every $i$. Also a straightforward double counting
of the pairs $(i,F)$, $i\in F$, $F\in \cF$ gives 
$dn\equiv kn \pmod p$. As $n\not \equiv 0$, we obtain 
$d\equiv k\pmod p$. Substituting this into (\ref{congr}), we infer 
\begin{equation} \label{congr2}
 k-\lambda\equiv -\lambda n + k^2 \pmod p,
\end{equation}
which, after rearrangement, gives the congruence to be proved. 
$\qed$

{\bf Proof of Theorem \ref{main5}:}
For each $i$ we define
$$
f_i(\ve x):= \ve x \cdot \ve v_i-\lambda \in \R[\ve x],
$$
where $\ve v_i\in \R^n$ is the characteristic vector of $F_i$. 
Let $P$ denote the space of all linear polynomials from $\R[\ve x]$. 
Then we have $\dim_{\R} P=n+1$.

It is easy to verify that the set of polynomials 
$\{f_i:~ 1\leq i\leq n\}\subset P$ 
is
linearly 
independent over $\R$ (substituting of $\ve v_i$ show that). Also, 
it is easily seen that the constant 
polynomial $1$ is not a linear combination of the polynomials $f_i$. 
One can readily verify this by considering a hypothetical relation 
$1=\sum_{i=1}^n\alpha_if_i$, and substituting 
$\ve v_i$, and $\ve 0 =(0,\ldots,0)$. 

Thus the set of polynomials $B:=\{1\}\cup \{f_i:~ 1\leq i\leq n\}$ is 
linearly 
independent over $\R$. Since $|\cF|=n$ and $\dim_{\R} P=n+1$, 
we see that $B$ is a basis of $P$.
From this we infer also that the homogeneous linear polynomials 
$\ve x \cdot \ve v_i$, $i=1, \ldots ,n$  are linearly independent over
$\R$. Let $A$ stand for the $n$ by $n$ matrix whose rows are the vectors
$\ve v_1,\ldots ,\ve v_n$. By the preceding observation $A$ is a nonsingular
matrix. 

We expand the monomials $x_i$ in the basis 
$B$:
\begin{equation} \label{exp}
x_i=\sum_{j=1}^n \theta_{i,j} f_j+ \kappa_i
\end{equation} 
for each $i$. 
Substituting $\ve v_j$ into (\ref{exp}) we find that
$$
\theta_{i,j}=(1-\kappa_i)/(|F_j|-\lambda)
$$
if $i\in F_j$, and
$$
\theta_{i,j}=(-\kappa_i)/(|F_j|-\lambda)
$$
if $i\not\in F_j$. These give 
\begin{equation} \label{exp2}
x_i=(1-\kappa_i) \sum_{j:i\in F_j} \frac{f_j}{|F_j|-\lambda}-\kappa_i
\sum_{j:i\not\in F_j} \frac{f_j}{|F_j|-\lambda} +\kappa_i.
\end{equation} 

If we compare the coefficients of $x_i$ on the two sides of (\ref{exp2}), 
we obtain $1>\kappa_i$ and 
\begin{equation} \label{exp5}
\sum_{j:i\in F_j} \frac{1}{|F_j|-\lambda}=\frac{1}{1-\kappa_i}.
\end{equation}

Let $r_i$ stand for the number of indices $j$ such that $i \in F_j$.
Substitute $\ve 1 = (1,\ldots ,1)$ into (\ref{exp2}). We obtain 
$$ 1=(1-\kappa_i)r_i-\kappa_i(n-r_i)+\kappa_i$$ 
and in turn 
\begin{equation} \label{fokos}
r_i=\kappa_i(n-1)+1~ \mbox{for}~ i=1,\ldots ,n. 
\end{equation}


Substitute the vector $\ve 0=(0,\ldots ,0)$ into (\ref{exp2}). Then we have
\begin{equation} \label{exp3}
0=(1-\kappa_i)(-\lambda) \sum_{j:i\in F_j} \frac{1}{|F_j|-\lambda}-\kappa_i  
(-\lambda)\sum_{j:i\not\in F_j} \frac{1}{|F_j|-\lambda} +\kappa_i.
\end{equation} 
Using (\ref{exp3}) and (\ref{exp5}) it is
easy to verify that $\kappa_i\ne 0$.
It follows then from (\ref{exp3}) and (\ref{exp5}), that 
\begin{equation} \label{exp6}
\sum_{j:i\not\in F_j} \frac{1}{|F_j|-\lambda}=\frac{1}{\kappa_i}-\frac{1}{\lambda}.
\end{equation} 
From (\ref{exp5}) and (\ref{exp6}) we arrive to 
\begin{equation} \label{exp7}
\sum_{j=1}^n  \frac{1}{|F_j|-\lambda}=\frac{1}{\kappa_i} +\frac{1}{1-\kappa_i} -\frac{1}{\lambda}.
\end{equation}
But the sum on the left side does not depend on $i$, implying that 
the numbers $\kappa_i$ all satisfy the same quadratic equation 
which follows 
from (\ref{exp7}). Let the roots of this equation be $\kappa $ and $\kappa'$. 
First suppose that $\kappa_i=\kappa$ for each $i$. 
Then $r:=r_i=\kappa (n-1)+1$
for each $i$. Let $y_1,\ldots ,y_n$ be a set of variables and consider the
following system of linear equations: 
\begin{equation} 
\sum _{j:~i\in F_j}y_j= \frac{1}{1-\kappa}~~~~~~ \mbox{for}~ i=1,\ldots, n. 
\end{equation}
The matrix of the system is the nonsingular $A^T$, hence the system is 
uniquely solvable. We readily see (using that the numbers $r_i$ are of the same
value $r$) that $y_j= \frac{1}{(1-\kappa)r}$ is a solution ($j=1,\ldots,
n$). But (\ref{exp5}) shows that  $y_j=\frac{1}{|F_j|-\lambda}$ 
is also a solution. 
We conclude that the sets $F_j$ all have the same size, and then this size
is necessarily $r$. This gives alternative (A) of the theorem. 

As for alternative (B), 
suppose that $\kappa_i$ takes on  the two different values $\kappa$ and
$\kappa'$, as $i$ runs over $[n]$. Then for the corresponding degrees $r$
and $r'$ we have by (\ref{fokos}) and $\kappa+\kappa'=1$ the following
$$ r+r'=\kappa(n-1)+1+\kappa'(n-1)+1=n+1, $$
proving the theorem. $\Box$

\section{A concluding remark}

It would be interesting to find a relation among the parameters of a
uniform extremal family $\cF\subset 2^{[n]}$ of $w$ element sets whose
elements have two different intersections $l_1$ and $l_2$. We note that
the conjectured size of such $\cF$ is ${n-w+2 \choose 2}$ in Conjecture 1.3
of \cite{BGKLTY}. See also the related bounds in Proposition 2.6,
Corollary 4.2, Theorem 5.5, and Proposition 6.2 of \cite{BGKLTY}.


\end{document}